\documentclass[12pt]{article}

\usepackage{amssymb,amsmath,graphics}

\numberwithin{equation}{section}

\allowdisplaybreaks

\newenvironment{Proof}{\removelastskip\par\medskip
\noindent{\em Proof.}
\rm}{\penalty-20\null\hfill$\square$\par\medbreak}

\newtheorem{prop}{Proposition}
\newtheorem{lemma}{Lemma}
\newtheorem{definition}{Definition}

\newtheorem{remark}{Remark}

\numberwithin{prop}{section} \numberwithin{lemma}{section}
\numberwithin{corollary}{section} \numberwithin{definition}{section}
\numberwithin{theorem}{section} \numberwithin{remark}{section}

\begin{document}

 ~
 \begin{center}

{\Large A stochastic volatility model with jumps}

\baselineskip=0.7cm
\end{center}

\begin{center}
{\large Youssef El-Khatib\footnote{%
Youssef\_Elkhatib@uaeu.ac.ae} \ } {\normalsize \ }

{\normalsize \ }
\emph{\ United Arab Emirates University\\ Department of Mathematical
Sciences, P.O.Box 17551, Al-Ain, U.A.E\ }
\end{center}
{\normalsize \ }
\begin{abstract}
\baselineskip0.5cm We consider a stochastic volatility model with jumps
where the underlying asset price is driven by the process sum of a
2-dimensional Brownian motion and a 2-dimensional compensated Poisson process.
The market is incomplete, resulting in infinitely many Equivalent Martingale
Measures. We find the set equivalent martingale measures, and we hedge by minimizing the variance 
using Malliavin calculus.
\end{abstract}

~

\baselineskip=0.5cm \noindent\textbf{Keywords:} stochastic volatility model, jumps, European options, incomplete markets, Malliavin calculus, mean-variance hedging.

\baselineskip0.7cm

\section{Introduction}

In the pioneer work of Black and Scholes (1973), the
financial asset prices are modeled by the Brownian motion, which
is a continuous process. The black and Scholes model does
not take into account the jumps which can occur at any time and randomly.
Three years later, Merton (1976) suggested a model with jumps. Since then,
the study of financial mathematical models have attracted the
interest of many mathematicians.

\noindent The recent international financial crisis and
its effects on global stock markets have showed once again the importance
of adding jumps to financial modeling for stock prices. Unlike the continuous case, discontinuous models assume this powerful
hypothesis: at any moment, a financial price can jump to decrease (increase) and attain,
in a negligible time, a significant lower (higher) value. In other words,
these models can simulate financial crisis, thus their importance.

\noindent On the other hand, the models in Black and Scholes (1973) and Merton (1976) assume a
deterministic volatility. Later on, new models with stochastic
volatility have been suggested to take into account the
so called \emph{smile} effect. Most of the works on these models assume
-for simplification- the continuity of the asset price trajectories
(driven by a Brownian motion). For continuous models with stochastic volatility,
We refer the reader to\footnote{the list is not exhaustive.}
 Heston (1993), Hull and  White (1987), Stein and Stein (1991) and Hagan, Kumar, Lesniewski and Woodward (2002).

\noindent We need more realistic models where the stochastic process describing the price trajectories involves
jumps. And in the same models, volatility should be stochastic not deterministic in order to consider the \emph{smile} effect. Several papers on stochastic volatility models including jumps
have been done. These works show clearly that the stochastic volatility models combined with
jump-diffusions are the best for modeling stock prices. Nevertheless they are still not well explored. This is due to their complication. Actually, they generate the incompleteness of the market, i.e., not every contingent claim can be hedged. For instance, while in Bates (1996) the stock price
dynamics includes jumps, the stochastic volatility is still considered continuous. In
Duffie, Pan, and Singleton (2000) or Broadie, Chernov, and Johannes (2005) the stochastic volatility contains jumps. However, these papers do not deal with the problem of equivalent martingale measure nor with the problem of hedging strategies for options.

\noindent In this work we are interested in a more general framework for discontinuous dynamics for the asset price with discontinuous stochastic volatility. The main contribution of this work is solving two problems: finding the equivalent martingale measure minimizing the entropy and finding hedging strategies under a general framework for jump-diffusions markets combined with stochastic volatility.

\noindent Assume that we have a market with two assets: a risky asset $S$ which is
related to a European Call option and a risk free one with price's process $%
A:=(A_{t})_{t\in \lbrack 0,T]}$ where $dA_{t}=r_{t}A_{t}dt,\ \ \ t\in
\lbrack 0,T],\ \ \ A_{0}>0$, and $r_{t}$ is a deterministic function
denoting the interest rate. Formally, let the underlying asset price of $S$
be
\begin{eqnarray*}
\frac{dS_{t}}{S_{t}} &=&\mu _{t}dt+\sigma
(t,Y_{t})[a_{t}^{(1)}dW_{t}^{(1)}+a_{t}^{(3)}dM_{t}^{(1)}],\ \ \ t\in
\lbrack 0,T],\ \ \ S_{0}=x>0, \\
&&\mbox{with} \\
dY_{t} &=&\mu _{t}^{Y}dt+\sum_{i=1}^{2}\sigma
_{t}^{(i)}[a_{t}^{(i)}dW_{t}^{(i)}+a_{t}^{(i+2)}dM_{t}^{(i)}],\ \ \
Y_{0}=y\in {\mathord{\mathbb R}},
\end{eqnarray*}%
where $W=(W^{(1)},W^{(2)})$ is a 2-dimensional Brownian motion and $%
M=(M^{(1)},M^{(2)})$ is a 2-dimensional compensated Poisson process with
independent components and deterministic intensity $(\int_{0}^{t}\lambda
_{s}^{(1)}ds,\int_{0}^{t}\lambda _{s}^{(2)}ds)$.  We assume that for $1\leq i\leq
4,a^{(i)}:[0,T]\longrightarrow {\mathord{\mathbb R}}$ is a deterministic
function.\newline
The most serious problem in a stochastic volatility model is incompleteness.
These models involve the existence of infinitely many
equivalent martingale measures (E.M.M.) i.e probabilities equivalent to the
historical one under which the discounted prices are martingales. First we
characterize the set of E.M.M.. We show that a probability $Q$ equivalent to
the historical probability $P$ is specified by its Radon-Nikodym density
w.r.t $P$
\begin{eqnarray*}
\rho _{T} &=&\prod_{i=1}^{2}\exp \left( \int_{0}^{T}\beta
_{s}^{(i)}dW_{s}^{(i)}-\frac{1}{2}\int_{0}^{T}(\beta
_{s}^{(i)})^{2}ds\right) \exp \left( \int_{0}^{T}\ln (1+\beta
_{s}^{(i+2)})dM_{s}^{(i)}\right.  \\
&&\left. +\int_{0}^{T}\lambda _{s}^{(i)}\left[ \ln (1+\beta
_{s}^{(i+2)})-\beta _{s}^{(i+2)}\right] ds\right) ,
\end{eqnarray*}%
where $(\beta _{t})_{t\in \lbrack 0,T]}$ is a ${\mathord{\mathbb R}}^{4}$%
-valued predictable process such that $\beta ^{(3)},\beta ^{(4)}>-1.$ If $Q$
is a $P-$E.M.M., $\beta ^{(1)}$ and $\beta ^{(3)}$ are related by
\begin{equation*}
\mu _{t}-r_{t}+\beta _{t}^{(1)}a_{t}^{(1)}\sigma (t,Y_{t})+\lambda
_{t}^{(1)}\beta _{t}^{(3)}a_{t}^{(3)}\sigma (t,Y_{t})=0,
\end{equation*}%
see Proposition~\ref{mmeq}.\newline
The process $\left( -\frac{\mu _{t}-r_{t}}{a_{t}^{(1)}\sigma (t,Y_{t})}%
,0,0,0\right) $ is an example of a ${\mathord{\mathbb R}}^{4}$-valued
predictable process satisfying the above equation, and it defines a $P$%
-E.M.M. This means that the set of $P$-E.M.M. is not empty. Moreover, since $%
\beta ^{(2)}$ and $\beta ^{(4)}$ do not appear in the last equation, so they
can be chosen arbitrarily, and thus there exists infinitely many $P-$E.M.M..

\subsection*{Mean-variance hedging}
We hedge using the mean-variance hedging approach initiated
by F\"{o}llmer and Sondermann (1986), %\cite{follmer4}
and we find the strategy by applying Malliavin calculus. \newline
Consider an option with payoff $f(S_T)$, where $(S_t)_{t\in [0,T]}$ is the
asset price with maturity $T$. We work with a $P$-E.M.M ${\hat Q}$. Let $(%
\hat{\eta}_t, \hat{\zeta}_t)_{t\in [0,T]}$ be a self-hedging strategy and $(%
\hat{V}_t)_{t\in [0,T]}$ be the portfolio value process. Using the chaotic
calculus, we conclude that the strategy minimizing the variance $E_{\hat{Q}}%
\left[(f(S_T)- \hat{V}_T )^2 \right]$ is given by
\begin{equation*}
{\hat{\eta}}_t =\frac{a^{(1)}_t E[D^{\hat{W}^{(1)}}_t f(S_T)\mid {\mathcal{F}%
}_t]+\lambda^{(1)}_t(1+{\hat{\beta}}^{(3)}_t)a^{(3)}_t E[D^{N^{(1)}}_t
f(S_T)\mid {\mathcal{F}}_t]} {((a^{(1)}_t)^2 +\lambda^{(1)}_t(1+{\hat{\beta}}%
^{(3)}_t)(a^{(3)}_t)^2)e^{\int_t^T r_sds}\sigma(t,Y_t) S_t},
\end{equation*}
where ${\hat{W}^{(1)}}_t={W}^{(1)}_t -\int_0^t \hat{\beta}^{(1)}ds$, and the
operators $D^{\hat{W}^{(1)}}$ and $D^{N^{(1)}}$ are respectively the
Malliavin derivative in the direction of the one dimensional Brownian motion
${\hat{W}}^{(1)}$ and the Malliavin operator in the direction of the Poisson
process $N^{(1)}$.\newline
This paper is organized as follows : In Section 2, we present some necessary
formulas. In the third section we introduce the model. The fourth one is
devoted to the hedging by minimizing the variance via Malliavin calculus. In
the last section, we characterize the E.M.M. minimizing the entropy, which
allows us to establish explicit formulae for the strategy.

\section{Preliminary}

Let $W=(W^{(1)},W^{(2)})$ be a 2-dimensional Brownian motion and $%
N=(N^{(1)},N^{(2)})$ be a 2-dimensional Poisson process with independent
components and deterministic intensity $(\int_0^t \lambda_s^{(1)}
ds,\int_0^t \lambda_s^{(2)} ds)$. We work in a probability space $(\Omega, {%
\mathcal{F}},({\mathcal{F}}_t)_{t\in [0,T]}, P)$, where $({\mathcal{F}}%
_t)_{t\in [0,T]}$ is the natural filtration generated by $W$ and $N$. We
denote by $M=(M^{(1)},M^{(2)})$ the associated compensated Poisson process,
i.e for $i=1,2$ and $t\in [0,T]$, we have $dM^{(i)}_t=dN^{(i)}_t
-\lambda^{(i)}_t dt$. Both $(\mathcal{F}_t)_{t\in [0,T]}$-martingales $W$
and $M$ are independent.

\begin{definition}
We denote by $\Gamma$ be the set of all $\mathcal{F}_t$-predictable
processes $(\gamma_t)_{t\in [0,T]}$ with values in ${\mathord{\mathbb R}}^4$
such that
\begin{equation*}
\sum_{i=1}^2 E_P \left[\int_0^t (\gamma^{(i)}_s)^2 ds\right] + \sum_{i=1}^2
E_P \left[\int_0^t (\gamma^{(i+2)}_s)^2 \lambda^{i}_s ds\right]<\infty, \ \
\ t\in [0,T].
\end{equation*}
\end{definition}

\noindent For a semi-martingale $X$ with $X_0=0$, the Dol\'eans-Dade exponential $%
\mathcal{E}(X)_t$ is the unique solution of the stochastic differential
equation
\begin{equation*}
Z_t=1+\int_0^t Z_{s^-} dX_s.
\end{equation*}
We have (Theorem~36 of Protter (1990))%
%\cite{protter5})
\begin{equation}  \label{expdade}
\mathcal{E}(X)_t=\exp\left(X_t -\frac{1}{2}[X_t,X_t]^c\right)\prod_{s\leq
t}(1+\Delta X_s)\exp(-\Delta X_s).
\end{equation}

\begin{remark}
Notice that for $\gamma \in \Gamma$ such that $\gamma^{(3)},\gamma^{(4)}>-1$
and for $i=1,2$
\begin{eqnarray*}
\mathcal{E}(\gamma^{(i)} W^{(i)})_t&=&\exp\left(\int_0^t \gamma^{(i)}_s
dW^{(i)}_s -\frac{1}{2}\int_0^t (\gamma^{(i)}_s)^2 ds\right), \\
\mathcal{E}(\gamma^{(i+2)} M^{(i)})_t&=&\exp\left(\int_0^t \ln
(1+\gamma^{(i+2)}_s) dM^{(i)}_s\right. \\
&&\left.+\int_0^t \lambda^{(i)}_s \left[\ln (1+\gamma^{(i+2)}_s)
-\gamma^{(i+2)}_s\right]ds\right).
\end{eqnarray*}
\end{remark}

\noindent The next lemma is the martingale representation theorem (Jacod (1979)).%
%\cite{jacod5}.

\begin{lemma}
\label{repres} Let $Z=(Z_t)_{t\in [0,T]}$ be a $\mathcal{F}_t$-martingale.
There exists a predictable process $\gamma \in \Gamma$ such that
\begin{equation*}
dZ_t = \sum_{i=1}^2 \gamma^{(i)}_t dW^{(i)}_t +\sum_{i=1}^2 \gamma^{(i+2)}_t
dM^{(i)}_t,\ \ \ t \in [0,T].
\end{equation*}
\end{lemma}

\section{The model}

\label{the model} Consider a market with two assets: a risky asset
which is related to a European call option and a riskless one. The maturity is $%
T$ and the strike is $K$. The price of the riskless asset is given by
\begin{equation*}
dA_t=r_t A_t dt, \ \ \ t\in [0,T],\ \ \ A_0=1,
\end{equation*}
where $r_t$ is deterministic and denotes the interest rate. The price of the
risky asset has a stochastic volatility and is given by
\begin{eqnarray*}
\frac{dS_t}{S_t}&=&\mu_t dt+\sigma(t,Y_t)[a^{(1)}_t dW^{(1)}_t +a^{(3)}_t
dM^{(1)}_t],\ \ \ t\in[0,T],\ \ \ S_0=x>0, \\
dY_t&=&\mu^Y_t dt+\sum_{i=1}^2 \sigma^{(i)}_t [a^{(i)}_t dW^{(i)}_t
+a^{(i+2)}_t dM^{(i)}_t ],\ \ \ t\in[0,T],\ \ \ Y_0=y \in {\mathord{\mathbb
R}},  \notag
\end{eqnarray*}
where for $1\leq i\leq 4, a^{(i)}: [0,T] \longrightarrow {\mathord{\mathbb R}%
}$ is a deterministic function. We assume that
\begin{equation*}
\sigma(t,.)\neq 0,\ \ \ \mbox{and}\ \ \ 1+\sigma(t,.)a^{(3)}_t >0, \ \ \
t\in [0,T].
\end{equation*}
We have
\begin{eqnarray*}
S_{t} &=& x\exp \left( \int_0^t a^{(1)}_s \sigma(s,Y_s)dW^{(1)}_s + \int_0^t
(\mu_{s}- a^{(3)}_s \lambda^{(1)}_s \sigma(s,Y_s) - \frac{1}{2}
(a^{(1)}_s)^2 \sigma^2(s,Y_s)) ds \right) \\
&&\times \prod_{k=1}^{k=N_t} (1+a^{(3)}_{T^{(1)}_k}
\sigma(T^{(1)}_k,Y_{T^{(1)}_k})),
\end{eqnarray*}
$0\leq t\leq T$, where $(T^{(1)}_k)_{k\geq 1}$ denotes the jump times of $%
(N^{(1)}_t)_{t\in [0,T]}$.

\subsection{Change of probability}

Let $Q$ be a $P$-equivalent probability; by the Radon-Nikodym theorem there
exists a ${\mathcal{F}}_T$-measurable random variable, $\rho_T :=\frac{dQ}{dP%
}$, such that $Q(A)=E_P[\rho_T 1_A]$, $A \in \mathcal{P}(\Omega)$. Notice
that $\rho_T$ is strictly positive $P$-a.s, since $Q$ is equivalent to $P$,
and $E_P[\rho_T]=E_P[\rho_T 1_\Omega]=1$. Consider now the $P$-martingale $%
\rho=(\rho_t)_{t\in [0,T]}$ defined by
\begin{equation*}
\rho_t :=E_P[\rho_T \mid \mathcal{F}_t]=E_P\left[\frac{dQ}{dP} \mid \mathcal{%
F}_t\right].
\end{equation*}

\begin{definition}
$\mathcal{H}$ is the set of all $P-$E.M.M., i.e $Q \in \mathcal{H}$ if and only
if $Q \simeq P$ and the discounted prices are $Q$-martingales.
\end{definition}

\noindent The next proposition gives the Radon-Nikodym density w.r.t $P$ of a $P$-E.M.M..

\begin{prop}
\label{mmeq}Let $Q \in \mathcal{H}$. There exists a predictable process $%
(\beta_t)_{t\in [0,T]}$ taking values in ${\mathord{\mathbb R}}^4$ such that
$\beta^{(3)},\beta^{(4)}>-1$ and the Radon-Nikodym density of $Q$ w.r.t $P$
is given by
\begin{eqnarray}  \label{rhob}
\rho_T&=&\prod_{i=1}^2 \mathcal{E}(\beta^{(i)}W^{(i)})_T \mathcal{E}%
(\beta^{(i+2)}M^{(i)})_T  \notag \\
&=& \prod_{i=1}^2 \exp\left(\int_0^T \beta^{(i)}_s dW^{(i)}_s -\frac{1}{2}%
\int_0^T (\beta^{(i)}_s)^2 ds\right)\exp\left(\int_0^T \ln
(1+\beta^{(i+2)}_s) dM^{(i)}_s\right.  \notag \\
&&\left.+\int_0^T \lambda^{(i)}_s \left[\ln (1+\beta^{(i+2)}_s)
-\beta^{(i+2)}_s\right]ds\right).
\end{eqnarray}
Moreover $\beta^{(1)}$ and $\beta^{(3)}$ are related by
\begin{equation}  \label{rel}
\mu_t -r_t +\beta^{(1)}_t a^{(1)}_t \sigma(t,Y_t) +\lambda^{(1)}_t
\beta^{(3)}_t a^{(3)}_t \sigma(t,Y_t)=0.
\end{equation}
\end{prop}

\begin{Proof}
We follow Bellamy (1999) %\cite{bellamy5}
for the case of a discontinuous market with deterministic volatility. By the
martingale representation theorem (Lemma~\ref{repres}) there exists a
predictable process $(\gamma_t)_{t \in [0,T]} \in \Gamma$ such that
\begin{equation*}
d\rho_t= \sum_{i=1}^2 \gamma^{(i)}_t dW^{(i)}_t +\sum_{i=1}^2
\gamma^{(i+2)}_t dM^{(i)}_t ,\ \ \ t \in [0,T].
\end{equation*}
We have $P(\rho_t >0, t\in [0,T])=1$; assuming $\beta:=\frac{\gamma}{\rho}$,
we obtain
\begin{eqnarray*}
\frac{d\rho_t}{\rho_t}&=&\sum_{i=1}^2 \beta^{(i)}_t dW^{(i)}_t
+\sum_{i=1}^2\beta^{(i+2)}_t dM^{(i)}_t,\ \ \ t \in [0,T].
\end{eqnarray*}
(\ref{rhob}) follows from (\ref{expdade}). In addition $(e^{-\int_0^t r_s
ds}S_t)_{t\in [0,T]}$ is a $Q$-martingale, in other words, $(e^{-\int_0^t
r_s ds}S_t \rho_t)_{t\in [0,T]}$ is a $P$-martingale. The integration by
parts formula (Protter (1990)) %\cite{protter5}
gives
\begin{equation*}
d(e^{-\int_0^t r_s ds}S_t \rho_t)=\rho_t d(e^{-\int_0^t r_s ds}S_t)
+e^{-\int_0^t r_s ds}S_t d\rho_t +d[e^{-\int_0^t r_s ds}S_t,\rho_t],
\end{equation*}
with
\begin{eqnarray*}
d[e^{-\int_0^t r_s ds}S_t,\rho_t]&=&\beta^{(1)}_t a^{(1)}_t \sigma(t,Y_t)dt
+ \beta^{(3)}_t a^{(3)}_t \sigma(t,Y_t)dN^{(1)}_t, \\
&=&\left(\beta^{(1)}_t a^{(1)}_t \sigma(t,Y_t)+\lambda^{(1)}_t \beta^{(3)}_t
a^{(3)}_t \sigma(t,Y_t)\right)dt+ \beta^{(3)}_t a^{(3)}_t
\sigma(t,Y_t)dM^{(1)}_t.
\end{eqnarray*}
Therefore
\begin{eqnarray*}
d(e^{-\int_0^t r_s ds}S_t \rho_t)&=&\rho_t S_t e^{-\int_0^t r_s ds}\left[%
(\mu_t -r_t +\beta^{(1)}_t a^{(1)}_t \sigma(t,Y_t) +\lambda^{(1)}_t
\beta^{(3)}_t a^{(3)}_t \sigma(t,Y_t))dt\right. \\
&&+\left.(\beta^{(1)}_t+\sigma(t,Y_t)a^{(1)}_t)dW^{(1)}_t +\beta^{(2)}_t
dW^{(2)}_t \right. \\
&&+\left.\left(\sigma(t,Y_t)a^{(3)}_t + \beta^{(3)}_t(1+
\sigma(t,Y_t)a^{(3)}_t)\right)dM^{(1)}_t+\beta^{(4)}_t dM^{(2)}_t\right].
\end{eqnarray*}
Thus $Q$ is a $P$-E.M.M. if
\begin{equation*}
\mu_t -r_t +\beta^{(1)}_t a^{(1)}_t \sigma(t,Y_t) +\lambda^{(1)}_t
\beta^{(3)}_t a^{(3)}_t \sigma(t,Y_t)=0.
\end{equation*}
\end{Proof}

\noindent Notice that there are no restrictions on $\beta^{(2)}$ and $\beta^{(4)}$,
which means that if $\mathcal{H} \neq \emptyset$, then $\mathcal{H}$
contains infinitely many $P$-E.M.M..

\section{Equivalent Martingale Measure minimizing the entropy}

Let $\Gamma^{\mathcal{H}}$ be the set of processes $\beta \in \Gamma$
satisfying (\ref{rel}). The Radon-Nikodym derivative $\rho_T$ associated to $%
\beta$ and given by (\ref{rhob}) defines a $P$-E.M.M.. From now on, a $P$-E.M.M. $%
Q $ in $\mathcal{H}$ will be denoted by $Q^{\beta}$, where $\beta \in
\Gamma^{\mathcal{H}}$. The process $\left(\frac{\mu_t-r_t}{a^{(1)}_t
\sigma(t,Y_t)},0,0,0 \right)$ belongs to $\Gamma^{\mathcal{H}}$ and it
defines a $P$-E.M.M., so $\mathcal{H}\neq \emptyset$. Thus $\mathcal{H}$
contains infinitely many $P$-E.M.M.. We choose the one that minimizes the
relative entropy. Let $Q^{\beta} \in \mathcal{H}$. Denoting by $%
I(Q^{\beta},P)$ the relative entropy of $Q^{\beta}$ w.r.t $P$, we have
\begin{equation*}
I(Q^{\beta},P)=E_P\left[\frac{dQ^{\beta}}{dP}\ln\frac{dQ^{\beta}}{dP}\right].
\end{equation*}
Our aim is to minimize $I(P,Q^{\beta})$ under $\mathcal{H}$. We have
\begin{equation*}
I(P,Q^{\beta})=E_{Q^{\beta}}\left[\frac{dP}{dQ^{\beta}}\ln\frac{dP}{%
dQ^{\beta}}\right]
\end{equation*}
Therefore the problem is to find a $\hat{\beta}$ which satisfies
\begin{equation}  \label{min}
I(P,Q^{\hat{\beta}}) = \min_{\beta \in \Gamma^{\mathcal{H}}}-E_P\left[\ln%
\frac{dQ^{\beta}}{dP}\right].
\end{equation}

\begin{lemma}
\label{lem} The minimization problem (\ref{min}) is equivalent to the
minimization of
\begin{eqnarray*}
&&(\mu_t -r_t+\lambda^{(1)}_t a^{(3)}_t \beta^{(3)}_t \sigma (t,Y_t))^2
-2\sigma^2 (t,Y_t) (a^{(1)}_t)^2 \lambda^{(1)}_t \left[\ln(1+\beta^{(3)}_t)-%
\beta^{(3)}_t \right] \\
&&-2\sigma^2 (t,Y_t) (a^{(1)}_t)^2 \lambda^{(2)}_t \left[\ln(1+%
\beta^{(4)}_t)-\beta^{(4)}_t\right],
\end{eqnarray*}
under all $\beta=\left(\frac{\mu_t -r_t+\lambda^{(1)}_t a^{(3)}_t
\sigma(t,Y_t)\beta^{(3)}_t }{\sigma(t,Y_t)a^{(1)}_t},0,\beta^{(3)},
\beta^{(4)}\right) \in \Gamma^{\mathcal{H}}.$
\end{lemma}

\begin{Proof}
Let $Q^{\beta} \in \mathcal{H}$. By (\ref{rhob})
\begin{eqnarray*}
I(P,Q^{\beta})&=&-E_P\left[\ln\frac{dQ^{\beta}}{dP}\right] \\
&=&E_P\left[\int_0^T \sum_{i=1}^2 \frac{1}{2}(\beta^{(i)}_t)^2
-\lambda^{(i)}_t \left[\ln(1+\beta^{(i+2)}_t)-\beta^{(i+2)}_t\right]dt\right]
\\
&=&E_P\left[\int_0^T \frac{G(\beta_t)}{2\sigma^2 (t,Y_t)(a^{(1)}_t)^2}dt%
\right],
\end{eqnarray*}
where $\beta \in \Gamma^{\mathcal{H}}$, and $G$ is the function defined by
\begin{eqnarray*}
G(\beta_t)&=& 2\sigma^2 (t,Y_t)(a^{(1)}_t)^2\left(\frac{1}{2}%
(\beta^{(1)}_t)^2 +\frac{1}{2}(\beta^{(2)}_t)^2 -\lambda^{(1)}_t \left[%
\ln(1+\beta^{(3)}_t)-\beta^{(3)}_t \right]\right. \\
&&\left.-\lambda^{(2)}_t \left[\ln(1+\beta^{(4)}_t)-\beta^{(4)}_t\right]%
\right),\ \ \ t\in [0,T].
\end{eqnarray*}
For a fixed $t$, we have by (\ref{rel}),
\begin{eqnarray*}
G(\beta_t)&=&(\mu_t -r_t+\lambda^{(1)}_t a^{(3)}_t
\sigma(t,Y_t)\beta_t^{(3)} )^2 +\sigma^2 (t,Y_t)
(a^{(1)}_t)^2\left((\beta^{(2)}_t)^2 \right. \\
&&\left.-2 \lambda^{(1)}_t \left[\ln(1+\beta^{(3)}_t)-\beta^{(3)}_t \right]%
-2 \lambda^{(2)}_t \left[\ln(1+\beta^{(4)}_t)-\beta^{(4)}_t\right]\right),\
\ \ t\in [0,T].
\end{eqnarray*}
Since $\beta^{(2)}_t$ appears only in the term $\sigma^2 (t,Y_t)
(a^{(1)}_t)^2(\beta^{(2)}_t)^2$ which is always positive, $\beta^{(2)}_t$
must be equal to zero.
\end{Proof}

\noindent The following proposition gives the solution to the minimization \ref{min}.

\begin{prop}
\label{mmeqmva} Consider $({\hat{\beta}}^{(1)}_t,{\hat{\beta}}^{(2)}_t, {%
\hat{\beta}}^{(3)}_t, {\hat{\beta}}^{(4)}_t)_{t \in[0,T]} \in \Gamma^{%
\mathcal{H}}$, with
\begin{equation*}
{\hat{\beta}}^{(2)}_t={\hat{\beta}}^{(4)}_t=0,\ \ \ {\hat{\beta}}^{(1)}_t=
\begin{cases}
\frac{r_t -\mu_t -\lambda^{(1)}_t a^{(3)}_t \sigma(t,Y_t) {\hat{\beta}}%
^{(3)}_t}{\sigma(t,Y_t) a^{(1)}_t} & \ \ \ \mbox{if} \ \ \ a^{(1)}\neq 0, \\
0 & \ \ \ \mbox{if} \ \ \ a^{(1)}=0,%
\end{cases}%
\end{equation*}
and let ${\hat{\beta}}^{(3)}_t$ be the unique solution of the equation
\begin{equation}  \label{entropie}
\lambda^{(1)}_t \sigma(t,Y_t)(a^{(3)}_t)^2 x+(a^{(1)}_t)^2
\sigma(t,Y_t)\left(\frac{x}{1+x}\right)-a^{(3)}_t (r_t-\mu_t)=0.
\end{equation}%
Then, the $P$-E.M.M. $\hat{Q}$ defined by its Radon-Nikodym density
\begin{equation*}
\prod_{i=1}^2 \mathcal{E}({\hat{\beta}}^{(i)}W^{(i)})_T \mathcal{E}({\hat{%
\beta}}^{(i+2)}M^{(i)})_T,
\end{equation*}
is the $P$-E.M.M. minimizing $I(P,Q^{\beta})$.
\end{prop}
\begin{Proof}
By Lemma~\ref{lem}, we have to minimize the function $F :]-1,\infty[\times]%
-1,\infty[ {\longrightarrow {\mathord{\mathbb R}}}$ defined by
\begin{eqnarray*}
F(x,y)&=&(\mu_t -r_t+\lambda^{(1)}_t a^{(3)}_t \sigma(t,Y_t)x)^2-2\sigma^2
(t,Y_t) (a^{(1)}_t)^2 \left(\lambda^{(1)}_t \left[\ln(1+x)-x\right]\right. \\
&&+\left. \lambda^{(2)}_t \left[\ln(1+y)-y\right]\right),
\end{eqnarray*}
for a fixed $t$ in $[0,T]$. Let $F^{^{\prime }}_{x}$ and $F^{\prime }_y$
denote the first order partial derivatives of $F$. The critical points of $F$
are determined by solving the equations $F^{^{\prime }}_{x} (x,y) =
F^{\prime }_y (x,y) = 0.$ Let $\hat{x}$ be the solution of (\ref%
{entropie}). It is unique since the function
\begin{equation*}
x \longrightarrow 2(\lambda^{(1)}_t)^2 \sigma^2(t,Y_t)(a^{(3)}_t)^2
x+2(a^{(1)}_t)^2 \lambda^{(1)}_t \sigma^2(t,Y_t)\frac{x}{1+x}%
+2\lambda^{(1)}_t \sigma(t,Y_t) a^{(3)}_t (\mu_t -r_t),
\end{equation*}
is strictly increasing from $]-1,\infty[$ to ${\mathord{\mathbb R}}$. One
can check that $({\hat{x}},0)$ is the only point which satisfies $F^{\prime
}_x (x,y) = F^{\prime }_y (x,y) = 0.$ Moreover we have
\begin{equation*}
(F^{^{\prime \prime }}_{xy}(\hat{x},0))^2 -F^{^{\prime \prime }}_{x^2}(\hat{x%
},0) F^{^{\prime \prime }}_{y^2}(\hat{x},0)<0\ \ \ \mbox{and}\ \ \
F^{^{\prime \prime }}_{x^2}(\hat{x},0)>0.
\end{equation*}
Therefore $F$ has a strict local minimum at $(\hat{x},0)$. This minimum is
global since $F$ goes to infinity when $x$ ($y$) approaches infinity.
\end{Proof}

\section{Hedging}

\label{hed} In this section we are interested in finding an optimal hedging strategy
for the model described in Section~\ref{the model}. We find the
strategy minimizing the variance using the Malliavin calculus. From now on,
we work with $\hat{Q}$: the $P$-E.M.M. minimizing the entropy given by
$\hat{\beta}$ from Proposition.~\ref{mmeqmva}. Consider
the two processes $\hat{W}=(\hat{W}^{(1)},\hat{W}^{(2)})$ and $\hat{M}=(\hat{%
M}^{(1)},\hat{M}^{(2)})$ where for $i=1,2$
\begin{equation*}
{\hat{W}^{(i)}}_t={W}^{(i)}_t -\int_0^t {\hat{\beta}}^{(i)}ds,\ \ \ t\in
[0,T],\ \ \ \mbox{and} \ \ \ {\hat{M}^{(i)}}_t={M}^{(i)}_t -\int_0^t
\lambda^{(i)}_s {\hat{\beta}}^{(i+2)}_s ds,\ \ \ t\in [0,T],
\end{equation*}
by Girsanov theorem (Jacod (1979)) %\cite{jacod5}
$\hat{W}$ is a ${\hat{Q}}$-Brownian motion and $\hat{M}$ is a ${\hat{Q}}$%
-compensated Poisson process. Under ${\hat{Q}}$, $(S_t)_{t\in[0,T]}$
satisfies
\begin{equation*}
\frac{dS_t}{S_t}=r_t dt+\sigma(t,Y_t)[a^{(1)}_t d{\hat{W}^{(1)}}_t
+a^{(3)}_t d{\hat{M}^{(1)}}_t],\ \ \ t\in[0,T],\ \ \ S_0=x>0.
\end{equation*}

\subsection{Chaotic calculus}

Let us denote by $\hat{X}$ the process
\begin{equation*}
(\hat{X}_{t}^{(1)},\hat{X}_{t}^{(2)},\hat{X}_{t}^{(3)},\hat{X}_{t}^{(4)})=(\hat{W}_{t}^{(1)},\hat{W}_{t}^{(2)},\hat{M}_{t}^{(1)},\hat{M}_{t}^{(2)}),\ \ \ t\in \lbrack 0,T],
\end{equation*}%
and let $({\hat{\mathcal{F}}}_{t})_{t\in \lbrack 0,T]}$ be the natural
filtration generated by $\hat{X}$.
\noindent We define the multiple stochastic integral and introduce the Malliavin
gradient and the Clark-Ocone formula in the multidimensional
Brownian-Poisson case (the following definitions and formulas can be
extended for the $d-$dimensional case, $d\leq 4$). For more details we refer
to L{\o }kka (1999), Nualart (1995), Nualart and Vives(1990), {\O }ksendal
(1996) and Privault (1997 a,b).
%\cite{lokka2},\cite{Nual2}, \cite{nualart2}, \cite{oksendal}, \cite{pri97a} or \cite{pri97b}.
Let $(e_{1},e_{2},e_{3},e_{4})$ be the canonical base of ${%
\mathord{\mathbb
R}}^{4}$. For $g_{n}\in L^{2}([0,T]^{n})$ we define the $n$-th iterated
stochastic integral of the function $f_{n}e_{i_{1}}\otimes \ldots \otimes
e_{i_{n}}$, with $1\leq i_{1},\ldots ,i_{n}\leq 4$, by
\begin{equation*}
I_{n}(g_{n}e_{i_{1}}\otimes \ldots \otimes
e_{i_{n}}):=n!\int_{0}^{T}\int_{0}^{t_{n}}\ldots
\int_{0}^{t_{2}}g_{n}(t_{1},\ldots ,t_{n})d{\hat{X}}_{t_{1}}^{(i_{1})}\ldots
d{\hat{X}}_{t_{n}}^{(i_{n})}.
\end{equation*}%
The iterated stochastic integral of a symmetric function $%
f_{n}=(f_{n}^{(i_{1},\ldots ,i_{n})})_{1\leq i_{1},\ldots ,i_{n}\leq 4}\in $
\newline
$L^{2}([0,T],{\mathord{\mathbb R}}^{4})^{\otimes n}$, where $%
f_{n}^{(i_{1},\ldots ,i_{n})}\in L^{2}([0,T]^{n})$, is
\begin{eqnarray*}
I_{n}(f_{n}):= &&\sum_{i_{1},\ldots ,i_{n}=1}^{4}I_{n}(f_{n}^{(i_{1},\ldots
,i_{n})}e_{i_{1}}\otimes \ldots \otimes e_{i_{n}}) \\
&=&n!\sum_{i_{1},\ldots ,i_{n}=1}^{4}\int_{0}^{T}\int_{0}^{t_{n}}\ldots
\int_{0}^{t_{2}}f_{n}^{(i_{1},\ldots ,i_{n})}(t_{1},\ldots ,t_{n})d{\hat{X}}%
_{t_{1}}^{(i_{1})}\ldots d{\hat{X}}_{t_{n}}^{(i_{n})}.
\end{eqnarray*}%
Recall that $\hat{X}$ has the Chaotic Representation Property (CRP) which
states that any square-integrable ${\hat{%
\mathcal{F}}}_{T}$-measurable functional can be expanded into a series of
multiple stochastic integrals -w.r.t ${\hat{X}}_{t}$- of deterministic
functions. For $F\in L^{2}(\Omega )$, there exists a unique sequence $(f_{n})_{n\in {%
\mathord{\mathbb N}}}$ of deterministic symmetric functions $%
f_{n}=(f_{n}^{(i_{1},\ldots ,i_{n})})_{i_{1},\ldots ,i_{n}\in \{1,\ldots
,4\}}\in L^{2}([0,T],{\mathord{\mathbb R}}^{4})^{\circ n}$ such that
\begin{equation}
F=\sum_{n=0}^{\infty }I_{n}(f_{n}).  \label{prc}
\end{equation}
\begin{definition}
Let $l\in \{1,\ldots ,4\}$, the operator ${\hat{D}}^{(l)}:{\mathrm{\mathrm{%
Dom~}}}({\hat{D}}^{(l)})\subset L^{2}(\Omega )\rightarrow L^{2}(\Omega
,[0,T])$ maps $F\in {\mathrm{\mathrm{Dom~}}}({\hat{D}}^{(l)})$ ($F$ having
the representation (\ref{prc})) to the process $({\hat{D}}_{t}^{(l)}F)_{t\in
\lbrack 0,T]}$ given by
\begin{eqnarray*}
\lefteqn{\hat{D}_{t}^{(l)}F:=\sum_{n=1}^{\infty
}\sum_{h=1}^{n}\sum_{i_{1},\ldots ,i_{n}=1}^{4}1_{\{i_{h}=l\}}} \\
&&I_{n-1}(f_{n}^{(i_{1},\ldots ,i_{n})}(t_{1},\dots ,t_{l-1},t,t_{l+1}\ldots
,t_{n})e_{i_{1}}\otimes \ldots \otimes e_{i_{h-1}}\otimes e_{i_{h+1}}\ldots
\otimes e_{i_{n}}) \\
&=&\sum_{n=1}^{\infty }nI_{n-1}(f_{n}^{l}(\ast ,t)),\ \ \ dP\times dt-a.e.
\end{eqnarray*}%
with $f_{n}^{l}=(f_{n}^{(i_{1},\ldots ,i_{n-1},l)}e_{i_{1}}\otimes \ldots
\otimes e_{i_{n-1}})_{1\leq i_{1},\ldots ,i_{n-1}\leq {4}}$.
\end{definition}
\noindent The domain of ${\hat{D}}^{(l)}$ is
\begin{eqnarray*}
{\mathrm{\mathrm{Dom ~}}}({\hat{D}}^{(l)})&=&\left\{F=\sum_{n=0}^{\infty}%
\sum_{i_1,\ldots,i_n=1}^{4} I_n(f_n^{(i_1,\ldots,i_n)}e_{i_1}\otimes \ldots
\otimes e_{i_n}) \in L^2(\Omega):\right. \\
&&\left.\sum_{i_1,\ldots,i_n=1}^{4} \sum_{n=0}^{\infty}n
n!\|f_n^{(i_1,\ldots,i_n)}\|^2_{L^2([0,T]^n)}< \infty\right\}.
\end{eqnarray*}
The probabilistic interpretations of ${\hat{D}}^{(l)}$ for
the Brownian motion and Poisson process cases are respectively given below.
\begin{description}
\item[The Brownian operator] For $1\leq l\leq 2$, the operator ${\hat{D}}%
^{(l)}$ is, in fact, the Malliavin derivative in the direction of the one
dimensional Brownian motion ${\hat{W}}^{(l)}$. So, we have for $1\leq l\leq
2 $ and $F=f\left( {\hat{W}}_{t_{1}},\ldots ,{\hat{W}}_{t_{n}}\right) \in
L^{2}(\Omega )$, where $(t_{1},\ldots ,t_{n})\in \lbrack 0,T]^{n}$ and $%
f(x^{11},x^{21},\ldots ,x^{1n},x^{2n})\in \mathcal{C}_{b}^{\infty }({%
\mathord{\mathbb R}}^{2n})$
\begin{equation*}
{\hat{D}}_{t}^{(l)}F=\sum_{k=1}^{k=n}\frac{\partial f}{\partial x^{lk}}%
\left( {\hat{W}}_{t_{1}},\ldots ,{\hat{W}}_{t_{n}}\right) 1_{[0,t_{k}]}(t).
\end{equation*}%
To find the Mallaivin derivative of an It\^{o} integral, we need the
following proposition (see corollary~5.13 of {\O }ksendal (1996)).%
% \cite{oksendal}).

\begin{prop}
\label{derivint} Let $(u_{t})_{t\in \lbrack 0,T]}$ be a ${\mathcal{\hat{F}}}%
_{t}-$adapted process such that $u_{t}\in {\mathrm{\mathrm{Dom~}}}({\hat{D}}%
^{(l)})$. Then for $l=1,2$ we have
\begin{equation*}
{\hat{D}}_{t}^{(l)}\int_{0}^{T}u_{s}d{\hat{W}}_{s}^{(l)}=\int_{t}^{T}({\hat{D%
}}_{t}^{(l)}u_{s})d{\hat{W}}_{s}^{(l)}+u_{t},
\end{equation*}
\end{prop}

\item[The Poisson operator] For $3\leq l\leq 4$, ${\hat{D}}^{(l)}$ is the
Malliavin operator\footnote{%
Notice that, unlike the Brownian case, the Malliavin operator in the Poisson
space is not a derivative.} in the direction of the Poisson process $%
N^{(l-2)}$. For $F\in {\mathrm{\mathrm{Dom ~}}}({\hat{D}}^{(l)})$
\begin{equation*}
{\hat{D}}^{(l)}_t F(\omega^{(1)},\ldots,\omega^{4})=\left\{
\begin{array}{ll}
F(\omega^{(1)},\omega^{(2)},\omega^{(3)}+1_{[t,\infty[},\omega^{(4)}
)-F(\omega^{(1)},\ldots,\omega^{(4)}), & l=3, \\
F(\omega^{(1)},\omega^{(2)},\omega^{(3)},\omega^{(4)}+1_{[t,\infty[}
)-F(\omega^{(1)},\ldots,\omega^{(4)}), & l=4.%
\end{array}
\right.
\end{equation*}
\end{description}

\noindent The Clark-Ocone formula is given by the next proposition.

\begin{prop}
\textbf{(The Clark-Ocone formula)} \label{proCO} Consider a
square-integrable, ${\hat{\mathcal{F}}}_{T}$-measurable, functional $F$ such
that $F\in \bigcap_{l=1}^{4}{\mathrm{\mathrm{Dom~}}}(\hat{D}^{(l)})$. $F$
has the following predictable representation
\begin{equation*}
F=E_{\hat{Q}}[F]+\sum_{l=1}^{2}\int_{0}^{T}E_{{\hat{Q}}}[{\hat{D}}%
_{t}^{(l)}F\mid {\hat{{\mathcal{F}}}}_{t}]d{\hat{W}}_{t}^{(l)}+%
\sum_{l=1}^{2}\int_{0}^{T}E_{{\hat{Q}}}[{\hat{D}}_{t}^{(l+2)}F\mid {\hat{{%
\mathcal{F}}}}_{t}]d{\hat{M}}_{t}^{(l)}.
\end{equation*}
\end{prop}

\subsection{Strategy minimizing the variance}

Suppose that we are required to find a portfolio $(\hat{\zeta}_{t},\hat{\eta}_{t})_{t\in \lbrack
0,T]}$ which leads to a given value $\hat{V}_{T}=F$. The process $(\hat{V}_{t})_{t\in \lbrack
0,T]}$ denote the value of the portfolio and $\hat{\zeta}_{t}$ and $\hat{\eta}_{t}$ denote the number of shares invested at time $t$ in the risky and in the riskfree assets respectively. We have for, $t\in
\lbrack 0,T]$, $\hat{V}_{t}=\hat{\zeta}_{t}A_{t}+\hat{\eta}_{t}S_{t}$. The strategy is assumed to be self-financing thus $d\hat{V}_{t}=\hat{\zeta}%
_{t}dA_{t}+\hat{\eta}_{t}dS_{t}$ and
\begin{equation*}
d\hat{V}_{t}=r_{t}\hat{V}_{t}dt+\sigma (t,Y_{t})\hat{\eta}%
_{t}S_{t}[a_{t}^{(1)}d{\hat{W}^{(1)}}_{t}+a_{t}^{(3)}d{\hat{M}^{(1)}}%
_{t}],\quad t\in \lbrack 0,T].
\end{equation*}%
Moreover for any $t\in \lbrack 0,T]$, we have
\begin{eqnarray*}
d\left( e^{\left( -\int_{0}^{t}r_{s}ds\right) }\hat{V}_{t}\right)
&=&-r_{t}e^{\left( -\int_{0}^{t}r_{s}ds\right) }\hat{V}_{t}dt+e^{\left(
-\int_{0}^{t}r_{s}ds\right) }d\hat{V}_{t} \\
&=&e^{\left( -\int_{0}^{t}r_{s}ds\right) }\left[ -r_{t}\hat{V}_{t}dt+r_{t}%
\hat{V}_{t}dt+\sigma (t,Y_{t})\hat{\eta}_{t}S_{t}[a_{t}^{(1)}d{\hat{W}^{(1)}}%
_{t}+a_{t}^{(3)}d{\hat{M}^{(1)}}_{t}]\right] ,
\end{eqnarray*}%
therefore
\begin{equation*}
e^{\left( -\int_{0}^{T}r_{s}ds\right) }\hat{V}_{T}=\hat{V}%
_{0}+\int_{0}^{T}e^{\left( -\int_{0}^{t}r_{s}ds\right) }\sigma (t,Y_{t})\hat{%
\eta}_{t}S_{t}[a_{t}^{(1)}d{\hat{W}^{(1)}}_{t}+a_{t}^{(3)}d{\hat{M}^{(1)}}%
_{t}].
\end{equation*}%
and
\begin{equation}
\label{proval}
\hat{V}_{T}=\hat{V}_{0}e^{\left(\int_{0}^{T}r_{s}ds\right) }+\int_{0}^{T}e^{\left( -\int_{t}^{T}r_{s}ds\right)
}\sigma (t,Y_{t})\hat{\eta}_{t}S_{t}[a_{t}^{(1)}d{\hat{W}^{(1)}}%
_{t}+a_{t}^{(3)}d{\hat{M}^{(1)}}_{t}].
\end{equation}%
Assuming that $F$ satisfies the hypothesis of the Proposition~\ref%
{proCO}, apply the Clark-Ocone formula to $F$. Comparing with the
equation (\ref{proval}), we see that the equality $\hat{V}_{T}=F$ cannot
hold unless
\begin{equation}
 \label{cond}
E[D_{t}^{\hat{W}^{(2)}}F\mid {\hat{{\mathcal{F}}}}_{t}]=E[D_{t}^{\hat{N}%
^{(2)}}F\mid {\hat{{\mathcal{F}}}}_{t}]=0,
\end{equation}%
because the expression of $\hat{V}_{T}$ in (\ref{proval}) does not contain
an integral term w.r.t. $d{\hat{W}^{(2)}}_{t}$ nor w.r.t. $d{\hat{M}%
^{(2)}}_{t}$.\newline
Take $F$ equals to the payoff $f(S_{T})$ of the model in section.~\ref{the model}), we see that
(\ref{cond})is not satisfied, because $D_{t}^{\hat{W}^{(2)}}f(S_{T})=f^{^{\prime }}(S_{T})D_{t}^{\hat{W}%
^{(2)}}S_{T}\neq 0$ and $D_{t}^{\hat{N}^{(2)}}f(S_{T})\neq 0$. In other words, the payoff $f(S_{T})$ is not attainable.
The market is then incomplete.\newline
Next we aim to find the strategy $(\hat{\zeta}_{t},\hat{\eta}_{t})_{t\in \lbrack 0,T]}$
that minimizes the variance
\begin{equation}
E_{\hat{Q}}\left[ (f(S_{T})-\hat{V}_{T})^{2}\right] .  \label{vari}
\end{equation}%
The next proposition gives the strategy minimizing the variance for our model considered in the Section~\ref%
{the model}.
\begin{prop}
\label{malliavinstra} The strategy minimizing (\ref{vari}) in the model of
Section~\ref{the model} is given by
\begin{equation}
{\hat{\eta}}_{t}=\frac{a_{t}^{(1)}E[D_{t}^{\hat{W}^{(1)}}f(S_{T})\mid {\hat{{%
\mathcal{F}}}}_{t}]+\lambda _{t}^{(1)}(1+{\hat{\beta}}%
_{t}^{(3)})a_{t}^{(3)}E[D_{t}^{N^{(1)}}f(S_{T})\mid {\hat{{\mathcal{F}}}}%
_{t}]}{((a_{t}^{(1)})^{2}+\lambda _{t}^{(1)}(1+{\hat{\beta}}%
_{t}^{(3)})(a_{t}^{(3)})^{2})e^{\int_{t}^{T}r_{s}ds}\sigma (t,Y_{t})S_{t}}.
\label{couvmall}
\end{equation}
\end{prop}
\begin{Proof}
Notice that the payoff $f(S_{T})=(S_{T}-K)^{+}$ is ${\mathcal{\hat{F}}}_{T}$-measurable.
We approach the function $x\mapsto f(x)(=(x-K)^{+}\mbox{or}=(K-x)^{+})$
by polynomials on compact intervals and proceed as in {\O }ksendal (1996)pp.
5-13. %\cite{oksendal},
By dominated convergence, $f(S_{T})\in \bigcap_{l=1}^{4}{\mathrm{\mathrm{%
Dom~}}}(D^{(l)})$. Thus by applying the Clark-Ocone formula to $f(S_{T})$ and using (%
\ref{proval}) we obtain
\begin{eqnarray*}
\lefteqn{E_{\hat{Q}}\left[ f(S_{T})-\hat{V}_{T})^{2}\right] =} \\
&&E_{\hat{Q}}\left[ \left( \int_{0}^{T}\left( E[D_{t}^{\hat{W}%
^{(1)}}f(S_{T})\mid {\hat{{\mathcal{F}}}}_{t}]-e^{\int_{t}^{T}r_{s}ds}\sigma
(t,Y_{t})\hat{\eta}_{t}S_{t}a_{t}^{(1)}\right) d{\hat{W}^{(1)}}_{t}\right)
^{2}\right.  \\
&&+\left. \left( \int_{0}^{T}E_{\hat{Q}}[D_{t}^{\hat{W}^{(2)}}f(S_{T})\mid {%
\hat{{\mathcal{F}}}}_{t}]d{\hat{W}^{(2)}}_{t}\right) ^{2}+\left(
\int_{0}^{T}E_{\hat{Q}}[D_{t}^{{N}^{(2)}}f(S_{T})\mid {\hat{{\mathcal{F}}}}%
_{t}]d{\hat{M}^{(2)}}_{t}\right) ^{2}\right.  \\
&&+\left. \left( \int_{0}^{T}\left( E_{\hat{Q}}[D_{t}^{N^{(2)}}f(S_{T})\mid {%
\hat{{\mathcal{F}}}}_{t}]-e^{\int_{t}^{T}r_{s}ds}\sigma (t,Y_{t})\hat{\eta}%
_{t}S_{t}a_{t}^{(3)}\right) d{\hat{M}}^{1}(t)\right) ^{2}\right]  \\
&=&E_{\hat{Q}}\left[ \int_{0}^{T}h_{2}({\hat{\eta}}_{t})dt\right] ,
\end{eqnarray*}
where
\begin{eqnarray*}
h_{2}(x) &=&(E_{\hat{Q}}[D_{t}^{\hat{W}^{(2)}}f(S_{T})\mid {\hat{{\mathcal{F}%
}}}_{t}])^{2}+\lambda _{t}^{(2)}(1+{\hat{\beta}}%
_{t}^{(4)})(E[D_{t}^{N^{(2)}}f(S_{T})\mid {\hat{{\mathcal{F}}}}_{t}])^{2} \\
&&+\left( E_{\hat{Q}}[D_{t}^{\hat{W}^{(1)}}f(S_{T})\mid {\hat{{\mathcal{F}}}}%
_{t}]-e^{\int_{t}^{T}r_{s}ds}\sigma (t,Y_{t})xS_{t}a_{t}^{(1)}\right) ^{2} \\
&&+\lambda _{t}^{(1)}(1+{\hat{\beta}}_{t}^{(3)})\left( E_{\hat{Q}%
}[D_{t}^{N^{(1)}}f(S_{T})\mid {\hat{{\mathcal{F}}}}_{t}]-e^{%
\int_{t}^{T}r_{s}ds}\sigma (t,Y_{t})xS_{t}a_{t}^{(3)}\right) ^{2}.
\end{eqnarray*}
It is easily verified that $h_{2}$ is convex, hence its minimum is reached
at $h_{2}^{\prime }(x)=0.$ Therefore the strategy minimizing the variance
is given by (\ref{couvmall}).
\end{Proof}
\subsection{Explicit formulae}
In order to derive explicit formulas for the strategy obtained in Proposition~\ref%
{malliavinstra}, we consider the following two special cases of the model in Section~\ref%
{the model}.:  a continuous stochastic volatility model with Brownian motion and a pure jumps stochastic volatility model with Poisson process.
\subsubsection{Brownian case}

Assume that $a^{(1)}_t=a^{(2)}_t=1$ and $a^{(3)}_t=a^{(4)}_t=0$, so $%
(S_t)_{0\leq t \leq T}$ depends on Brownian information only. Under ${\hat Q}
$, $(S_t)_{0\leq t \leq T}$ is given by
\begin{equation*}
S_t=x\exp\left(\int_0^t \left(r_s - \frac{\sigma^2 (s,Y_s)}{2}\right) ds +
\int_0^t \sigma (s,Y_s)d{\hat{W}^{(1)}}_s \right),
\end{equation*}
with
\begin{equation*}
Y_t=y+\int_0^t \left(\mu^Y_s +\sigma^{(1)}_s \frac{r_s-\mu_s}{\sigma(s,Y_s)}%
\right)ds+ \int_0^t \sigma^{(1)}_s d{\hat{W}^{(1)}}_t +\int_0^t
\sigma^{(2)}_s dW^{(2)}_s.
\end{equation*}
In the following proposition we compute the Malliavin derivative of the
payoff $(S_T-K)^+$. We can replace the result in the formula (\ref{couvmall}%
), and obtain an explicit formula for the strategy.

\begin{prop}
We have
\begin{eqnarray}  \label{derivs}
D^{\hat{W}^{(1)}}_t(S_T-K)^+&=&1_{\{S_T > K \}}S_T
\left(\sigma(t,Y_t)+\int_t^T \frac{\partial \sigma}{\partial y}(s,Y_s)D^{%
\hat{W}^{(1)}}_s Y_t d\hat{W}^{(1)}_s \right.  \notag \\
&&-\left. \int_t^T \sigma(s,Y_s)\frac{\partial \sigma}{\partial y}(s,Y_s)D^{%
\hat{W}^{(1)}}_t Y_s ds\right)
\end{eqnarray}
where
\begin{equation}  \label{derivy}
D^{\hat{W}^{(1)}}_t Y_s=\sigma^{(1)}_t \exp\left(-\int_t^s \sigma^{(1)}_u
\frac{r_u-\mu_u}{\sigma^2(u,Y_u)}du\right)\ \ \ s\in[t,T].
\end{equation}
\end{prop}

\begin{Proof}
By the chain rule for $D^{\hat{W}^{(1)}}_t$ and thanks to Proposition~\ref%
{derivint} we obtain
\begin{eqnarray*}
\lefteqn{ D^{\hat{W}^{(1)}}_t(S_T-K)^+ =} \\
&&1_{\{S_T > K \}}S_T \left(D^{\hat{W}^{(1)}}_t \int_0^T \left(r_s - \frac{%
\sigma^2 (s,Y_s)}{2}\right) ds +D^{\hat{W}^{(1)}}_t \int_0^T \sigma(s,Y_s)d{%
\hat{W}^{(1)}}_s\right) \\
&=&1_{\{S_T > K \}}S_T \left(-\int_t^T D^{\hat{W}^{(1)}}_t \frac{\sigma^2
(s,Y_s)}{2}ds + \int_t^T D^{\hat{W}^{(1)}}_t \sigma(s,Y_s)d{\hat{W}^{(1)}}_s
+\sigma(t,Y_t)\right),
\end{eqnarray*}
which gives (\ref{derivs}). Concerning the other derivative, we have for $%
0\leq t\leq s\leq T$
\begin{eqnarray*}
D^{\hat{W}^{(1)}}_t Y_s &=& \int_t^s D^{\hat{W}^{(1)}}_t \left(\mu^Y_u
+\sigma^{(1)}_u \frac{r_u-\mu_u}{\sigma(u,Y_u)}\right)du+ \sigma^{(1)}_t \\
&=&\sigma^{(1)}_t -\int_t^s \sigma^{(1)}_u \frac{r_u-\mu_u}{\sigma^2(u,Y_u)}%
D^{\hat{W}^{(1)}}_t Y_u du,
\end{eqnarray*}
So for $t$ fixed in $[0,T]$, the Malliavin derivative of $Y_s$ for $s\in
[t,T]$, $(D^{\hat{W}^{(1)}}_t Y_s)_{s\in [t,T]}$, satisfies a stochastic
differential equation, whose solution is precisely (\ref{derivy}).
\end{Proof}

\subsubsection{The Poisson case}

Similarly, as in the Brownian case, we aim to compute the quantity $D^{%
\hat{M}^{(1)}}_t (S_T-K)^+$ and replace the result in the expression of the
strategy in order to obtain an explicit formula for the Poisson case.
Suppose that we are working in the
Poisson space with a $2$-dimensional Poisson
process. The underlying asset price $(S_t)_{0\leq t \leq T}$ depends on
the Poisson process only. Hence we assume that $a^{(3)}_t=a^{(4)}_t=1$ and $%
a^{(1)}_t=a^{(2)}_t=0$. Under ${\hat Q}$, the dynamics of $(S_t)_{0\leq t
\leq T}$ is given by
\begin{equation*}
S_{t} = x\exp \left(\int_0^t \left(\mu_{s}+ \frac{r_s -\mu_s}{\sigma(s,Y_s)}%
\ln(1+\sigma(s,Y_s))\right)ds +\int_0^t \ln(1+\sigma(s,Y_s))d\hat{M}^{(1)}_s
\right),
\end{equation*}
for $t\in[0,T]$. The process $(Y_t)_{t\in [0,T]}$ under ${\hat Q}$, has the
representation
\begin{equation*}
Y_t=y+\int_0^t \left(\mu^Y_s +\sigma^{(1)}_s \frac{r_s-\mu_s}{\sigma(s,Y_s)}%
\right)ds+ \int_0^t \sigma^{(1)}_s d{\hat{M}^{(1)}}_t +\int_0^t
\sigma^{(2)}_s dM^{(2)}_s.
\end{equation*}

\begin{prop}
\begin{eqnarray*}
\lefteqn{D^{\hat{M}^{(1)}}_t (S_T-K)^+ =}  \notag \\
&&-(S_T-K)^+ + \left(\exp \left\{\int_t^T \left[\mu_{s}+ \frac{r_s -\mu_s}{%
\sigma(s,Y_s +\sigma^{(1)}_t )}\ln(1+\sigma(s,Y_s +\sigma^{(1)}_t))\right]%
ds\right.\right.  \notag \\
&&+\left.\left.\int_t^T \ln(1+\sigma(s,Y_s +\sigma^{(1)}_t))d\hat{M}^{(1)}_s
\right\}\times S_t (1+\sigma(t,Y_t +\sigma^{(1)}_t))-K\right)^+
\end{eqnarray*}
\end{prop}

\begin{Proof}
Using the probabilistic interpretation of $D^{\hat{M}^{(1)}}_t$ given below,
we obtain
\begin{equation*}
D^{\hat{M}^{(1)}}_t (S_T-K)^+=(S_T(\omega +1_{[t,T]})-K)^{+}
-(S_T(\omega)-K)^+ .
\end{equation*}
But
\begin{eqnarray*}
\lefteqn{S_T(\omega +1_{[t,T]})=x\exp \left(\int_0^t \left[\mu_{s}+ \frac{%
r_s -\mu_s}{\sigma(s,Y_s(\omega +1_{[t,T]}))}\ln(1+\sigma(s,Y_s(%
\omega+1_{[t,T]})))\right]ds\right.} \\
&&\left.+\int_0^t \ln(1+\sigma(s,Y_s(\omega+1_{[t,T]})))d\hat{M}%
^{(1)}_s\right) \\
&&\times\exp \left(\int_t^T \left[\mu_{s}+ \frac{r_s -\mu_s}{%
\sigma(s,Y_s(\omega +1_{[t,T]}))}\ln(1+\sigma(s,Y_s(\omega+1_{[t,T]})))%
\right]ds\right. \\
&&\left.+\int_t^T \ln(1+\sigma(s,Y_s(\omega+1_{[t,T]})))d\hat{M}^{(1)}_s
\right)\times(1+\sigma(t,Y_t(\omega+1_{[t,T]}))),
\end{eqnarray*}
and
\begin{eqnarray*}
Y_t(\omega+1_{[t,T]})&=&Y_t +\sigma^{(1)}_t,\ \ \ t\in [0,T], \\
Y_s(\omega+1_{[t,T]})&=&%
\begin{cases}
Y_s & \ \ \ \mbox{if}\ \ \ s\in [0,t[, \\
Y_s +\sigma^{(1)}_t & \ \ \ \mbox{if}\ \ \ s\in [t,T].%
\end{cases}%
\end{eqnarray*}
The proof is complete.
\end{Proof}

\end{document}